\documentclass[12pt]{amsart}
\usepackage{pstricks}

\newtheorem{fact}{Fact}[section]
\newtheorem{lemma}[fact]{Lemma}
\newtheorem{theorem}[fact]{Theorem}
\newtheorem{definition}[fact]{Definition}
\newtheorem{example}[fact]{Example}

\newtheorem{proposition}[fact]{Proposition}
\newtheorem{corollary}[fact]{Corollary}

\DeclareMathOperator{\CC}{{\mathbb C}}
\DeclareMathOperator{\N}{\mathbb N}\DeclareMathOperator{\Z}{\mathbb
Z} \DeclareMathOperator{\C}{\mathcal C }
\DeclareMathOperator{\A}{\mathcal A}
\DeclareMathOperator{\T}{\mathcal T}
\DeclareMathOperator{\U}{\mathcal U}
\DeclareMathOperator{\n}{\mathfrak{n}_{--}}
\DeclareMathOperator{\np}{\mathfrak{n}_{+}}
\DeclareMathOperator{\g}{\mathfrak{g}}

\DeclareMathOperator{\h}{\mathfrak{h}}
\DeclareMathOperator{\res}{Res} \DeclareMathOperator{\sym}{Sym}
\DeclareMathOperator{\id}{id} \DeclareMathOperator{\asym}{ASym}
\DeclareMathOperator{\sgn}{sgn} \DeclareMathOperator{\dlog}{d \log}
 \DeclareMathOperator{\ep}{\varepsilon}
\DeclareMathOperator{\coeff}{coeff}
\DeclareMathOperator{\str}{String}

\def\li#1{\begin{picture}(30,10)\put(3,3){\line(1,0){24}}\put(12,-3){${}_{#1}$}\end{picture}}
\def\lii#1{\begin{picture}(30,10)\put(3,3){\line(1,0){24}}\put(8,-3){${}_{#1}$}\end{picture}}

\title{Poincar\'e-Birkhoff-Witt expansions of the canonical elliptic differential form}

\author{G. Felder}
\address{Departement Mathematik, ETH-Zentrum, 8092 Z\"urich,
Switzerland}
\email{felder@math.ethz.ch}

\author{R. Rim\'anyi}
\address{Department of Mathematics, University of North Carolina at Chapel Hill, USA}
\email{rimanyi@email.unc.edu}

\author{A. Varchenko}
\address{Department of Mathematics, University of North Carolina at Chapel Hill, USA}
\email{anv@email.unc.edu}

\begin{document}
\thanks{\noindent Supported by
        NSF grants DMS-0405723 (2nd author), DMS-0244579 (3rd author) 
AMS Subject classification 33C67}

\begin{abstract} We study the canonical $U(\n)$-valued elliptic differential form, whose projections
to different Kac-Moody algebras are key ingredients of the
hypergeometric integral solutions of elliptic KZ differential
equations and Bethe ansatz constructions. We explicitly determine
the coefficients of the projections in the simple Lie algebras $A_r,
B_r, C_r, D_r$ in a conveniently chosen Poincar\'e-Birkhoff-Witt
basis. As an application we give a new formula for eigenfunctions of
Hamiltonians of the Calogero-Moser model.
\end{abstract}

\maketitle

\section{Introduction}

Let $\g$ be a simple Lie algebra with Cartan decomposition
$\g=\n\oplus \h\oplus \np$. Let $V$ be the tensor product
$V_{\Lambda_1}\otimes \ldots\otimes V_{\Lambda_n}$ of highest weight
$\g$-modules. The $V$-valued hypergeometric solutions of the
elliptic Knizhnik-Zamolodchikov differential equations have the form
(\cite{fv},~cf.~\cite{sv})
\begin{equation}\label{integral}
I(z,\lambda,\tau)=\int_{\gamma(z,\lambda,\tau)} \Phi(t,z,\tau)\cdot
\Omega^{V}(t,z,\lambda,\tau).
\end{equation}
Here $t=(t_1,\ldots,t_N)$, $z=(z_1,\ldots,z_n)$, $\tau\in\CC$ with
Im $\tau>0$, $\lambda\in \h$, $\Phi$ is an explicit scalar-valued
master function, $\gamma$ is a suitable cycle in $t$-space, and
$\Omega^{V}$ is a $V$-valued differential $k$-form in $t$-space.

The same $\Phi$ and $\Omega^V$ have applications to the Bethe Ansatz
method. It is known \cite{asymp, etkir, etkir2, fv,fv2} that the
values of $\Omega^V$ at the critical points of $\Phi$ (with respect
to $t$) give eigenfunctions of the Hamiltonians of the quantum
elliptic Calogero-Moser model \cite{gibher} (in the $\g=sl_{r+1}$
case).

For every $V =V_{\Lambda_1}\otimes \ldots\otimes V_{\Lambda_n}$, the
$V$-valued differential form $\Omega^V$ can be constructed out of a
single $U(\n)$-valued differential form $\Theta^{\g}$, where $U(\ )$
denotes the universal enveloping algebra. The form $\Theta^{\g}$ is
called the canonical elliptic differential form. In applications it
is important to have convenient formulas for $\Theta^{\g}$, and this
is the goal of the present paper.

The form $\Theta^{\g}$ (depending also on a vector $k\in \N^r$) is
defined as the projection of the `canonical element' $\Omega_k$
according to the scheme
\[ \begin{matrix}
\Omega_k & \in & \A(k)\otimes U[k] & \to  & \A(k)\otimes U(\n)[k] && \\
& & \downarrow & & \downarrow & \\
&& W\otimes U[k] & \to & W \otimes U(\n)[k] & \ni & \Theta^{\g}_k.
\end{matrix}\]
Here $\A(k)$ is a homogeneous component of an Orlik-Solomon algebra;
$U[k]$ is its dual, and $W$ is a space of certain differential forms
written in terms of theta functions.

In Section 2 we define and study $\A(k)\otimes U[k]$, and its
canonical element $\Omega_k$. In Sections 3 and 4 we find
Poincar\'e-Birkhoff-Witt formulas for the projection of $\Omega_k$
in $\A(k)\otimes U(\n)[k]$. In Section 5 we define representations
of $\A(k)$ in spaces $W$ of suitable differential forms. In
Section~6 we use the PBW expansion of $\Theta^{sl_{r+1}}_k$ to find
a new formula for eigenfunctions of the Hamiltonian of the
Calogero-Moser model. In the Appendix we describe the image of the
 faithful representation $\A(k)\to W$.

\section{Cohomology of the complement of
the discriminantal arrangement. Combinatorial codes and operations.}

\subsection{} \label{start}
An {\em ordered $p$-forest} is a graph with no cycles, with $p$
edges, and a numbering of its edges by the numbers $1,2,\ldots,p$.
Let $\A_n^p$ be the complex vector space generated by the ordered
$p$-forests on the vertex set $\{z,t_1,\ldots,t_n\}$ (a set of
symbols now), modulo the following two kinds of relations:
\begin{description}
\item[R1]
$T_1=-T_2$ if $T_1$ and $T_2$ have the same underlying graph, and
the order of their edges differ by a transposition;
\item[R2]
\begin{equation}\label{triangle}
\begin{picture}(205,20)(0,15)
\put(15,0){$*$} \put(18,3){\line(-1,2){15}} \put(0,30){$*$}
\put(18,3){\line(1,2){15}}\put(30,30){$*$} \put(3,15){$a$}
\put(28,15){$b$}\put(48,15){$+$}

\put(85,0){$*$} \put(88,3){\line(-1,2){15}} \put(70,30){$*$}
\put(73,33){\line(1,0){30}}\put(100,30){$*$} \put(73,15){$b$}
\put(85,34){$a$} \put(118,15){$+$}

\put(155,0){$*$} \put(143,33){\line(1,0){30}} \put(140,30){$*$}
\put(158,3){\line(1,2){15}}\put(170,30){$*$} \put(155,34){$b$}
\put(168,15){$a$} \put(190,15){$=0$}
\end{picture}\qquad\qquad (a,b\in\{1,\ldots,p\}),\end{equation}
\vskip  .4 true cm \noindent that is, the sum of three $p$-forests
that locally (i.e. their subgraphs spanned by 3 vertices) differ as
above, but are otherwise identical, is 0.
\end{description}

Let $z$ be a complex number, and consider $\CC^{n}$ with coordinates
$t_1,\ldots,t_n$. The discriminantal arrangement $\C^n=\C^n(z)$ in
the space $\CC^{n}$ is the collection of hyperplanes $t_i=t_j$
($1\leq i<j\leq n$) and $t_i=z$ ($1\leq i\leq n$). Let $\U_n$ denote
the complement of $\C^n$, i.e. $\U_n=\CC^n\setminus \cup_{H\in \C^n}
H$.

\begin{proposition}\label{comp} The cohomology group
$H^p(\U_n;\CC)$ is isomorphic to $\A_n^p$.
\end{proposition}

\begin{proof} The definition of  $\A(k)$ is a combinatorial code
for the description of $H^p(\U_n;\CC)$ by Arnold in \cite{arnold}.
\end{proof}

Now let us fix $r\in \N$, and $r$ non-negative integers
$k=(k_1,k_2,\ldots,k_r)$ with $\sum k_i=|k|$. The set of coordinates
\begin{equation}\label{vars}
(t_1^{(1)},\ldots,t_{k_1}^{(1)},t_1^{(2)},\ldots,t_{k_2}^{(2)},
\ldots,t_1^{(r)},\ldots,t_{k_r}^{(r)}).
\end{equation}
will be denoted by $\T(k)$. We will consider $\CC^{|k|}$ with
coordinates $\T(k)$, and the space $\A_{|k|}^p$ with vertex set
$\T(k)\cup\{z\}$. The group $G_k=\prod \Sigma_{k_i}$ acts on
$\CC^{|k|}$ (by permuting the coordinates with the same upper
indices) which then induces an action of $G_k$ on $\A_{|k|}^p$. The
skew-invariant subspace (i.e. the collection of $x$'s for which
$\pi\cdot x=\sgn(\pi)x, \forall \pi\in G_k$) of $\A_{|k|}^{|k|}$
will be denoted by $\A(k)$.

\begin{example} \label{vecspace} \rm In the $r\leq 2$ examples (in the whole paper)
we will use $t_j$ for $t^{(1)}_j$ and $s_j$ for $t^{(2)}_j$. The
vector space $\A((1,1))$ is two-dimensional with a basis $\{
z\li{1}t_1\li{2}s_1, z\li{1}s_1\li{2}t_1\}$. The vector space
$\A((2))$ is one-dimensional with a basis $\{t_1\li{1}z\li{2}t_2\}$.
\end{example}

\subsection{The star-product.}

\begin{definition} Let $e_1,e_2,\ldots,e_{|k|}$ be the edges (in order)
of an ordered $|k|$-tree $T$ on the vertex set $\T(k)\cup\{z\}$. The
sign $\ep(T)$ of $T$ is defined by
$${\bigwedge_{i=1}^{|k|}d\ h(e_i)}=\ep(T)\cdot{\bigwedge_{i=1}^r\bigwedge_{j=1}^{k_i}d\ t^{(i)}_j},$$
where $h(e)$ is the `head' of the edge $e$, i.e. the end vertex of
$e$ farther from the root $z$.
\end{definition}

For example the signs of the three trees of Example \ref{vecspace}
are +1, -1, +1, respectively.

\begin{definition}
For an element $x$ in a vector space acted upon by $G_k$, let
$\asym_k(x)=\asym(x)$ be the anti-sym\-met\-ri\-za\-ti\-on of $x$:
$\sum_{\pi\in G_k} \sgn(\pi)\pi(x)$. Let $\sym_k(x)=\sym(x)$ be the
symmetrization of $x$: $\sum_{\pi\in G_k} \pi(x)$.
\end{definition}

For $k,l\in \N^r$ let $T_1$ and $T_2$ be ordered $|k|$- and
$|l|$-trees on the vertex set $\T(k)\cup\{z\}$ and $\T(l)\cup\{z\}$
respectively. We will define their {\em star}-product $T_1*T_2$ by
the following procedure. We replace the vertices $t^{(i)}_j$ of
$T_2$ with $t^{(i)}_{k_i+j}$, add $|k|$ to the numbers on the edges
of $T_2$, and identify the $z$-vertices of the two trees---thus we
obtain a tree $T$ with $|k+l|+1$ vertices. The product $T_1*T_2$ is
defined as $\ep(T_1)\ep(T_2)\ep(T)/\prod_i(k_i!l_i!)\cdot
\asym_{k+l}T$, i.e. $T_1* T_2=$

$$\frac{\ep(T_1)\ep(T_2)\ep(T)}{\prod_{i=1}^r
k_i!l_i!}\cdot\asym_{k+l}
\begin{pspicture}(7,6)(-3,2)
\pscurve{-}(2,1)(-2,2)(-1,4)(2,1) \pscurve{-}(2,1)(6,2)(5,4)(2,1)
\uput[dl](-0.5,6){$T_1$} \uput[dl](5.5,6){$T_2$}
\uput[dl](0,3){$t^{(i)}_j$} \uput[dl](5,3){$t^{(i)}_{k_i+j}$}
\uput[dl](1,5){(edges: $1,\ldots,|k|$)} \uput[dl](7,5){(edges:
$|k|+1,\ldots,|k+l|$)}\uput[dl](2.2,1.2){$\bullet$}
\uput[dl](2.3,0.6){$z$}
\end{pspicture}
$$

\vskip 2 true cm

The star-product induces a product (that we will also call
star-product) $\A(k)\otimes \A(l)\to \A(k+l)$. The vector space
$\oplus_{k\in \N^r} \A(k)$ with the star-product is an associative
and commutative algebra.

\begin{example} \rm
\[(z\li{1}t_1)*(z\li{1}s_1)=(t_1\li{1}z\li{2}s_1)=
(z\li{1}t_1\li{2}s_1)-(z\li{1}s_1\li{2}t_1).\] Computation shows
that
\[
(z\li{1}t_1)^{*n}=
n!\asym (z\li{1}t_{1}\li{2}t_{2}\li{3}\ldots\li{n}t_{n}).
\]
\end{example}

\subsection{Residue.}\label{residue} Let $T$ be a $p$-tree on the vertex set
$\T(k)\cup\{z\}$, and choose a coordinate $t^{(i)}_j$. We will
define the $t^{(i)}_j$-residue $\res^{(i)}_j T$ of the tree $T$ as a
$p-1$-tree on the vertex set $\T(k)\cup\{z\}\setminus\{t^{(i)}_j\}$,
constructed as follows.
\begin{itemize}
\item If $t^{(i)}_j$ and $z$ are connected by an edge with number
$a$, then we contract that edge (i.e. delete the edge, and identify
its vertices), and label this contracted vertex by $z$.  Then, if
the number of an edge is greater than $a$, we decrease it by 1.
$\res^{(i)}_j$ is defined as $(-1)^{a-1}$ times this modified tree.
See Figure 1.

\[ \res^{(i)}_j
\begin{pspicture}(7,4)(-1,1)
\pscurve{-}(2,0)(0,1.5)(0,3)(2,0) \pscurve{-}(4,0)(6,1.5)(6,3)(4,0)
\uput[dl](2.25,0.25){$\bullet$} \uput[dl](4.25,0.25){$\bullet$}
\put(2,0){\line(1,0){2}} \put(2,-.6){$z$} \put(4,-.6){$t^{(i)}_j$}
\put(3,.1){$a$}
\end{pspicture}=(-1)^{a-1}
\begin{pspicture}(5,4)(-.5,1)
\put(0,4){edges:}\put(1.6,4){$a+1\mapsto a$}
\put(1.6,3.5){$a+2\mapsto a+1$} \put(1.6,3){$\ldots$}
 \pscurve{-}(2,0)(0,1.5)(0,3)(2,0)
\pscurve{-}(2,0)(4,1.5)(4,3)(2,0) \uput[dl](2.25,0.25){$\bullet$}
\put(2,-.6){$z$}
\end{pspicture}\]

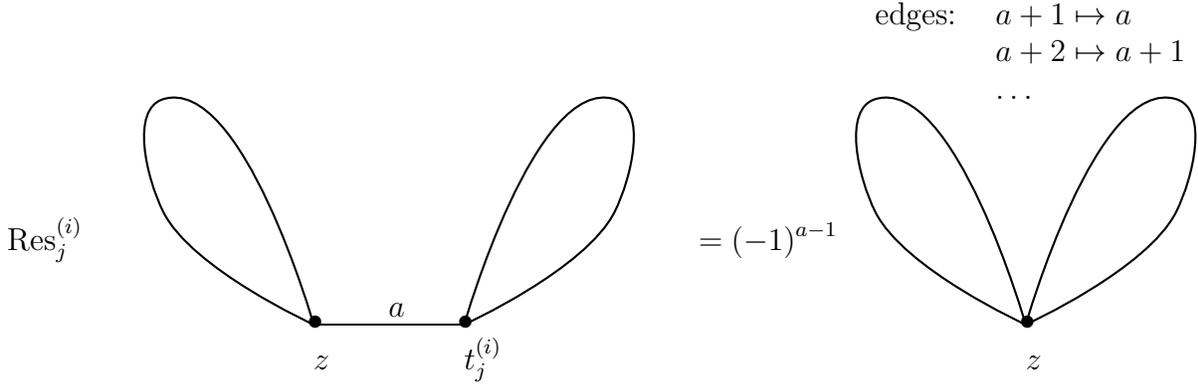
\vskip 1.5 true cm\centerline{Figure 1} \bigskip

\item If $t^{(i)}_j\li{}z$ is not an edge of $T$,  $\res^{(i)}_j T$ is defined to be $0$.
\end{itemize}

This operation is consistent with relations R1 and R2, hence they
are defined for $\A(k)$. Observe that $\res^{(i)}_{k_i}$ is a map
$\A(k)\to\A(k-1_i)$. (Here, and in the sequel, $1_i$ denotes a
vector, whose coordinates are 0 except the $i$'th coordinate, which
is 1.)

\subsection{The dual space of the cohomology of the complement.}

Let $U_r$ denote the free associative complex algebra generated by
the $r$ symbols $\tilde{f}_1,\ldots,\tilde{f}_r$. It is multigraded
by $\N^r$ ($\deg \tilde{f}_i=1$), the degree $k$ subspace will be
denoted by $U[k]$. In \cite{sv} $U[k]$ is shown to be isomorphic to
the $G_k$ skew-invariant part of the dual space of
$H^{|k|}(\U_{|k|};\CC)$. In this section we recall this isomorphism
in our combinatorial language.

\begin{definition}Let $r\in \N$ and $k\in \N^r$.
\begin{itemize}
\item Maps $J:\{1,\ldots,|k|\}\to\{1,\ldots,r\}$ with
$\#J^{-1}(i)=k_i$ will be called ($k$-)multiindices.

\item For a
multiindex $J$, let $\tilde{f}_J$ be the monomial
$\tilde{f}_{J(|k|)}\tilde{f}_{J(|k|-1)}\cdots \tilde{f}_{J(1)}$ in
$U[k]$.

\item For a multiindex $J$ let $c$ be the unique map
$\{1,\ldots,|k|\}\to \N$ that restricted to $J^{-1}(i)$ is the {\em
increasing} function onto $\{1,\ldots,k_i\}$.

\item The sign $\sgn(J)$ of a multiindex $J$ is $(-1)^m$ where $m$ is
the minimal number of transpositions to be applied to the list
$J(1),J(2),\ldots,J(|k|)$ to get
$1,\ldots,1,2,\ldots,2,\ldots,r,\ldots,r$.

\end{itemize}
\end{definition}

Let $T_0$ be the graph with a single vertex $z$ and with no edge.
For $T\in \A(k)$ and $J$ a $k$-multiindex, the expression
\begin{equation}\label{resdef}\sgn(J)\cdot
\res^{(J(|k|))}_{c(|k|)}\Big( \res^{(J(|k-1|))}_{c(|k-1|)}\big(
\ldots \res^{(J(1))}_{c(1)}(T)\ldots\big)\Big)
\end{equation}
is equal to a constant times $T_0$. Let this constant be denoted by
$<\tilde{f}_J,T>$.

\begin{theorem} \cite[Th.~2.4]{sv} \label{dualcomp}
The map $\tilde{f}_J\mapsto <\tilde{f}_J,\ \ \ >$ is an isomorphism
between $U[k]$ and the dual space $\A(k)^*$ of $\A(k)$.
\end{theorem}\qed

\begin{corollary} \label{coro}
For the basis $\{\tilde{f}_J\}$ of $U_r$ the
dual basis is
\[\tilde{f}_J^*=\sgn(J)\cdot \asym\Big(z \li{1}
t_{c(1)}^{(J(1))} \li{2} t_{c(2)}^{(J(2))} \li{3} \cdots \li{|k|}
t_{c(|k|)}^{(J(|k|))}\Big).\]
\end{corollary}

\subsection{Product on $\sum \A(k)$, coproduct on $\sum U[k]$.}

Recall that $U_r$ is equipped with its standard Hopf algebra
structure. The comultiplication $\Delta:U_r\to U_r\otimes U_r$ is
defined by $\Delta(x)=1\otimes x+x\otimes 1$ for degree 1 elements.

\begin{theorem} \label{duality}
Under the duality between $\sum_{k\in\N^r}\A(k)$ and
$U_r$ the coproduct $\Delta$ corresponds to the star-multiplication.
\end{theorem}

\begin{proof} We need the following notion: a triple
\[(S_1\subset\{1,\ldots,|k+l|\},S_2\subset\{1,\ldots,|k+l|\},J:\{1,\ldots,|k+l|\}\to\{1,\ldots,r\})\]
is called a {\em shuffle} of $J_1:\{1,\ldots,|k|\}\to
\{1,\ldots,r\}$ and $J_2:\{1,\ldots,|l|\}\to \{1,\ldots,r\}$ if
\begin{itemize}
\item $\#S_1=k$, $\#S_2=l$, $S_1\cup S_2=\{1,\ldots,|k+l|\}$,
\item for the increasing bijections $s_1:S_1\to \{1,\ldots,k\}$, $s_2:S_2\to \{1,\ldots,l\}$
we have
$$J(i)=\begin{cases} J_1\circ s_1(i) & i\in S_1 \\
J_2\circ s_2(i) & i\in S_2. \end{cases}$$
\end{itemize}

The definition of $\Delta$ implies that
$\Delta^*(\tilde{f}_{J_1}^*\otimes \tilde{f}_{J_2}^*)=\sum
\tilde{f}_J^*$, where the sum runs over all the shuffles of $J_1$
and $J_2$. On the other hand $\tilde{f}_{J_1}^* *
\tilde{f}_{J_1}^*=$
\[\pm\asym(z \li{1} t_{c_1(1)}^{(J_1(1))} \li{2} \cdots \li{|k|} t_{c_1(|k|)}^{(J_1(|k|))}) *
\asym(z \li{1} t_{c_2(1)}^{(J_2(1))} \li{2} \cdots \li{|k|}
t_{c_2(|k|)}^{(J_2(|k|))}).\] Using relation (\ref{triangle}), this
is equal to the sum of terms
$$\pm\asym(z \li{1} t_{c(1)}^{(J(1))} \li{2} \cdots \lii{|k+l|}
t_{c(|k+l|)}^{(J(|k+l|))}),
$$
where $J$ is a shuffle of $J_1$ and $J_2$. We illustrate this
argument with the following example:
$(\tilde{f}_2\tilde{f}_1)^**\tilde{f}_1^*=(z\li{1}t_1\li{2}s)*(z\li{1}t_1)=$
\[-\asym\Big(
\begin{picture}(65,30)
\put(5,0){$z$}\put(12,7){\line(1,1){10}}\put(13,14){${}_1$}
\put(28,17){$t_1$} \put(36,22){\line(1,0){15}}\put(42,25){${}_2$}
\put(53,17){$s$}
\put(12,-3){\line(1,-1){10}}\put(13,-14){${}_3$}\put(28,-20){$t_2$}
\end{picture}
\Big)=-\asym\Big(
\begin{picture}(65,30)
\put(5,0){$z$}\put(12,7){\line(1,1){10}}\put(13,14){${}_1$}
\put(28,17){$t_1$} \put(36,22){\line(1,0){15}}\put(42,25){${}_2$}
\put(53,17){$s$}
\put(35,11){\line(1,-1){10}}\put(36,2){${}_3$}\put(51,-4){$t_2$}
\end{picture}
+
\begin{picture}(90,30)
\put(5,0){$z$}\put(37,-10){\line(1,1){10}}\put(38,-3){${}_1$}
\put(53,0){$t_1$} \put(61,5){\line(1,0){15}}\put(67,8){${}_2$}
\put(78,0){$s$}
\put(12,-3){\line(1,-1){10}}\put(13,-14){${}_3$}\put(28,-20){$t_2$}
\end{picture}\Big)=
\]

\[
-\asym\Big(\begin{picture}(92,30)
\put(5,0){$z$}\put(12,7){\line(1,1){10}}\put(13,14){${}_1$}
\put(28,17){$t_1$} \put(36,22){\line(1,0){15}}\put(42,25){${}_2$}
\put(53,17){$s$}
\put(62,13){\line(1,-1){10}}\put(61,2){${}_3$}\put(76,-4){$t_2$}
\end{picture}
+
\begin{picture}(95,30)
\put(5,0){$z$}\put(12,7){\line(1,1){10}}\put(13,14){${}_1$}
\put(28,17){$t_1$} \put(61,0){\line(1,0){15}}\put(67,3){${}_2$}
\put(78,-5){$s$}
\put(35,11){\line(1,-1){10}}\put(36,2){${}_3$}\put(51,-4){$t_2$}
\end{picture}+
\begin{picture}(90,30)
\put(5,0){$z$}\put(37,-10){\line(1,1){10}}\put(38,-3){${}_1$}
\put(53,0){$t_1$} \put(61,5){\line(1,0){15}}\put(67,8){${}_2$}
\put(78,0){$s$}
\put(12,-3){\line(1,-1){10}}\put(13,-14){${}_3$}\put(28,-20){$t_2$}
\end{picture}\Big)=
\]
\ \vskip .8 true cm

\noindent $2\asym(z\li{1} t_1 \li{2} t_2 \li{3}
s)-\asym(z\li{1}t_1\li{2}s\li{3}t_2)=2(\tilde{f}_2\tilde{f}_1^2)^*+(\tilde{f}_1\tilde{f}_2\tilde{f}_1)^*$.

A careful examination of signs shows that $\Delta^*$ and the
star-product agree on the basis $\tilde{f}_J^*$, hence the Theorem
is proved. \end{proof}

\subsection{The canonical element.}
We obtained that the spaces $\A(k)$ and $U[k]$ are dual vector
spaces, and hence there is a canonical element in $\A(k)\otimes
U[k]$. It is defined by $\sum b_i^*\otimes b_i$ for a basis $b_i$ in
$U[k]$ and its dual basis $b_i^*$ in $\A(k)$. This definition does
not depend on the choice of the basis. Equivalently, the canonical
element is the identity in $\A(k)\otimes
U[k]=$Hom$(\A(k)^*,U[k])=$Hom$(U[k],U[k])$.

\begin{definition}
The canonical element in $\A(k)\otimes U[k]$ will be denoted by
$\Omega_k$.
\end{definition}

Using Corollary \ref{coro} we obtain an explicit form for
$\Omega_k$.

\begin{theorem} \label{canonical} We have
\begin{equation} \label{univ}
\Omega_k=\sum_{J} \sgn(J)\cdot \asym\Big(z \li{1} t_{c(1)}^{(J(1))}
\li{2} t_{c(2)}^{(J(2))} \li{3} \cdots \li{|k|}
t_{c(|k|)}^{(J(|k|))}\Big)\otimes \tilde{f}_J,
\end{equation}
where the summation runs over all $k$-multiindices.
\end{theorem}

\begin{example}
$\Omega_{(1,1)}=(z\li{1}t\li{2}s)\otimes
\tilde{f}_2\tilde{f}_1+(z\li{1}s\li{2}t)\otimes
\tilde{f}_2\tilde{f}_1$.
\begin{equation}
\begin{aligned}
\Omega_{(2,1)}= & \Big( (z\li{1}t_1\li{2}t_2\li{3}s)
-(z\li{1}t_2\li{2}t_1\li{3}s)\Big)\otimes \tilde{f}_2\tilde{f}_1^2- \\
&\Big( (z\li{1}t_1\li{2}s\li{3}t_2)
-(z\li{1}t_2\li{2}s\li{3}t_1)\Big)\otimes \tilde{f}_1\tilde{f}_2\tilde{f}_1+ \\
& \Big( (z\li{1}s\li{2}t_1\li{3}t_2)
-(z\li{1}s\li{2}t_2\li{3}t_1)\Big)\otimes \tilde{f}_2\tilde{f}_1^2.
\end{aligned}
\end{equation}
\end{example}

\section{Canonical elements associated with simple Lie algebras}

Let $\g$ be a simple Lie algebra of rank $r$ with Cartan
decomposition $\g=\n\oplus \h \oplus \np$. The $\N^r$-multigraded
universal enveloping algebra $U(\n)$ of $\n$ is generated by $r$
elements $f_1,\ldots,f_r$ (the standard Chevalley generators),
subject to the Serre relations. The degree $k$ part of $U(\n)$
($\deg f_i=1$) is denoted by $U(\n)[k]$. We have the quotient map
$q:U_r\to U(\n)$ by mapping $\tilde{f}_i$ to $f_i$.

\begin{definition}\label{qdef} The canonical element $\Omega^{\g}_k$ associated with the
simple Lie algebra $\g$ is defined as the image of the canonical
element $\Omega_k$ under the map
\begin{equation}
\id\otimes\ q: \A(k)\otimes U[k] \to \A(k)\otimes U(\n)[k].
\end{equation}
\end{definition}

\subsection{PBW expansion of the canonical element.}
The Lie algebra $\n$ is the direct sum of $1$-dimensional weight
spaces $\mathfrak{n}_\beta$ labeled by the positive roots:
$\n=\oplus_{\beta} \mathfrak{n}_{\beta}$. Now let us make the
following choices: \begin{itemize} \label{choices}
\item{} (C1) a linear ordering $\beta_1<\dots<\beta_m$ of the positive roots (thus
$m=\dim \n$),
\item (C2) a generator $F_{\beta_j}$ in
$\mathfrak{n}_{\beta_j}$. \end{itemize}
When considering
$F_{\beta_j}$ in $U(\n)$, let its degree be $k^{(j)}\in \N^r$.
Equivalently, $k^{(j)}_i$ is the coefficient of the $i$'th simple
root (i.e. the simple root corresponding to $f_i$) in the
decomposition of $\beta_j$. According to the
Poincar\'e-Birkhoff-Witt theorem, a $\CC$-basis of the algebra
$U(\n)$ is given by the collection of elements
$F_{\beta_1}^{p_1}F_{\beta_2}^{p_2}\ldots F_{\beta_m}^{p_m}$ where
$p=(p_1,\ldots,p_m)\in\N^m$.

Hence,  the canonical element associated with $\g$ can be written
(depending on the choices (C1,C2)) as
\begin{equation}
\Omega_k^{\g}\ =\sum_{p} T_p \otimes
  F_{\beta_1}^{p_1}F_{\beta_2}^{p_2}\ldots
F_{\beta_m}^{p_m},
\end{equation}
where the summation runs over those $p\in \N^m$ for which $\sum_j
p_jk^{(j)}=k\in \N^r$. Here $T_p\in \A(k)$ is a linear combination
of ordered spanning trees on the vertex set $\T(k)\cup\{z\}$.

\begin{example} \label{sl3} \rm In the Lie algebra $sl_3(\CC)$ of traceless
$3\times 3$ matrices let $f_1=\begin{pmatrix} 0 & 0 & 0 \\ 1 & 0 & 0
\\ 0 & 0 & 0 \end{pmatrix}$, $f_2=\begin{pmatrix} 0 & 0 & 0 \\ 0 & 0 & 0
\\ 0 & 1 & 0 \end{pmatrix}$. Let $\alpha_1$ and $\alpha_2$ be the simple
roots corresponding to $f_1$ and $f_2$ respectively. The positive
roots are $\alpha_1$, $\alpha_2$, and $\alpha_1+\alpha_2$. Let us
choose the ordering \begin{itemize} \item $\beta_1=\alpha_1\ \ <\ \
\beta_2=\alpha_1+\alpha_2\ \ <\ \ \beta_3=\alpha_2,$ \end{itemize}
\noindent and the elements \begin{itemize} \item $F_{\alpha_1}=f_1,
\ F_{\alpha_1+\alpha_2}=[f_2,f_1],\ F_{\alpha_2}=f_2,$
\end{itemize}
 with $k^{(1)}=(1,0)$, $k^{(2)}=(1,1)$, $k^{(3)}=(0,1)$. Then
for example $\Omega^{sl_3}_{(1,1)}$ can be written as
$T_{(1,0,1)}\otimes f_1f_2+T_{(0,1,0)}\otimes [f_2,f_1]$. Using
Theorem \ref{canonical} we obtain
\[\begin{aligned}
\Omega^{sl_3}_{(1,1)}= & (z\li{1}t\li{2}s)\otimes f_2f_1
-(z\li{1}s\li{2}t)\otimes f_1f_2=\\
& (z\li{1}t\li{2}s)\otimes
([f_2,f_1]+f_1f_2)-(z\li{1}s\li{2}t)\otimes f_1f_2=\\
& ((z\li{1}t\li{2}s)-(z\li{1}s\li{2}t))\otimes
f_1f_2+(z\li{1}t\li{2}s)\otimes [f_2,f_1].\end{aligned}\]
\end{example}

Observe that $T_{(1,0,1)}$ (the coefficient of $f_1f_2$) is equal to
$(z\li{1}t)*(z\li{1}s)$. In general we have the following formula.

\begin{theorem}{\rm (Product formula.)} \label{product}
 There exist unique elements
$\eta_{\beta_j}\in \A(k^{(j)})$ ($j=1,\ldots,m$) such that
\begin{equation}
T_p=\frac{1}{\prod_l
p_l!}\cdot\overbrace{\eta_{\beta_1}*\ldots*\eta_{\beta_1}}^{p_1}*
\overbrace{\eta_{\beta_2}*\ldots*\eta_{\beta_2}}^{p_2}* \ldots*
\overbrace{\eta_{\beta_m}*\ldots*\eta_{\beta_m}}^{p_m}.
\end{equation}
\end{theorem}

\begin{proof} Let $F^p=F_{\beta_1}^{p_1}F_{\beta_2}^{p_2}\ldots
F_{\beta_m}^{p_m}$. The key ob\-ser\-va\-tion is that in a PBW-basis
the co-mul\-ti\-pli\-ca\-tion is given by
\[ \Delta(F^p)=\sum_{p'+p''=p}\Big(\prod_{i=1}^m
\frac{p_i!}{p_i'!p_i''!}\Big) F^{p'}\otimes F^{p''}.\] Hence, for
its dual we get
\[\Delta^*(F^{p'*}\otimes F^{p''*})=\prod_{i=1}^m
\frac{p_i!}{p_i'!p_i''!}\cdot F^{(p'+p'')*}.\] Using Theorem
\ref{duality} this is equivalent to
$T_{p'}*T_{p''}=\prod_i(p_i!/p'_i!p''_i!)\cdot T_{p'+p''}$, from
which the theorem follows (put $\eta_{\beta_j}=T_{\beta_j}$).
\end{proof}

\begin{example} \rm
The computation in Example \ref{sl3} shows that (with the choices
made there) $\eta_{\alpha_1}=(z\li{1}t)$,
$\eta_{\alpha_1+\alpha_2}=(z\li{1}t\li{2}s)$,
$\eta_{\alpha_2}=(z\li{1}s)$. Then for example the coefficient of
$f_1^2[f_2,f_1]f_2$ in $\Omega^{sl_3}_{(3,2)}$ is
\[T_{(2,1,1)}=\frac{1}{2}\asym_{(3,2)}\Big(
\begin{picture}(80,30)(-10,20)
\put(0,20){$z$} \put(2,28){\line(0,1){15}} \put(0,45){$t_1$}
\put(5,26){\line(1,1){15}} \put(20,45){$t_2$}
\put(6,24){\line(1,0){15}} \put(24,20){$t_3$} \put(15,26){${}_3$}
\put(39,26){${}_4$} \put(33,24){\line(1,0){15}} \put(50,20){$s_1$}
\put(5,20){\line(1,-1){15}} \put(20,0){$s_2$}
\put(-3,34){${}_1$}\put(9,36){${}_2$} \put(8,8){${}_5$}
\end{picture}
\Big).
\]
\end{example}

\bigskip

\subsection{The residues of $\eta_{\beta}$.}

\begin{lemma} \label{fent} The map $\psi:U[k-1_i]\to U[k]$, $x\mapsto (-1)^{k_1+\ldots+k_i-1}x\cdot
\tilde{f}_i$ is dual to the map $\res^{(i)}_{k_i}:\A(k)\to
\A(k-1_i)$.
\end{lemma}

\begin{proof} Let $T\in\A(k)$, $f_J\in U[k-1_i]$. We need to check
that $<(-1)^{k_1+\ldots+k_i-1}\tilde{f}_J\tilde{f}_i,T>$
=$<\tilde{f}_J,\res^{(i)}_{k_i}T>$. Writing out definition
(\ref{resdef}) we obtain an identity.
\end{proof}

\begin{theorem} \label{resth} We have
\[\res^{(i)}_{k_i}\Omega_k=(-1)^{k_1+\ldots+k_i-1}\Omega_{k-1_i}\cdot(1\otimes
\tilde{f}_i).\]
\end{theorem}

\begin{proof} Let $\{b_u\}$ be a basis of $U_r[k-1_i]$, hence
$\Omega_{k-1_i}=\sum b_u^*\otimes b_u$. Since the map $\psi$ in
Lemma \ref{fent} is an embedding, the $\psi$ images of $b_u$ can be
extended to a basis $\{\psi(b_u), c_v\}$ of $U_r[k]$. Then
$\Omega_k= \sum \psi(b_u)^*\otimes \psi(b_u) + \sum c_v^*\otimes
c_v$. We have $(\res^{(i)}_{k_i}\otimes 1)\Omega_k=\sum
\res^{(i)}_{k_i}(\psi(b_u)^*) \otimes \psi(b_u)+\sum
\res^{(i)}_{k_i}(c_v^*)\otimes c_v$, which, according to the lemma,
is $\sum b_u^*\otimes \psi(b_u)=(1\otimes \psi)\Omega_{k-1_i}$, what
we wanted to prove.
\end{proof}

Recall that after making the choices (C1), (C2), any element in
$U(\n)$ can be written in the basis $\{F^p:p\in \N^m\}$ (recall
$F^p=F_{\beta_1}^{p_1}F_{\beta_2}^{p_2}\ldots F_{\beta_m}^{p_m}$).
If the element $x$ is expanded in this basis, let the coefficient of
$F^p$ be denoted by $\coeff(x,F^p)$. Theorem \ref{resth} has the
following direct corollary.

\begin{corollary} \label{coeffs} We have
\[(-1)^{k^{(j)}_1+\ldots+k^{(j)}_i-1}\res^{(i)}_{k^{(j)}_i} \eta_{\beta_j}=\sum_p
\coeff(F^p\cdot f_i,F_{\beta_j})\cdot T_p.\]
\end{corollary}

Since the residues of $\eta_{\beta}$'s will play a crucial role in
our argument below, we will make our choices so that the
coefficients on the right-hand side are easily computable.

\section{Computing $\eta_{\beta}$'s}

We reduced the computation of $T_p$'s to the computation of
$\eta_{\beta}$'s for the positive roots $\beta$. In this section we
compute the elements $\eta_{\beta}$ (with suitable choices (C1),
(C2)) for the simple Lie algebras of type A, B, C, D.

\begin{theorem} \label{barr} Let $\g$ be a simple Lie algebra of rank $r$,
with positive roots $\beta_1, \ldots, \beta_m$ and let choices (C1),
(C2) be made. Suppose $\bar{\eta}_{\beta_j}\in\A(k^{(j)})$
($j=1,\ldots,m$)  satisfy
\begin{description}
\item[(i)] $<\tilde{F},\bar{\eta}_{\beta_j}>=0$ if $\tilde{F}\in U[k^{(j)}]$ belongs to
the ideal of Serre relations;
\item[(ii)]
$\res^{(i)}_{k_i}\bar{\eta}_{\beta_j}=\res^{(i)}_{k_i}{\eta}_{\beta_j}$
for all $i=1,\ldots,r$.
\end{description}
Then $\bar{\eta}_{\beta_j}=\eta_{\beta_j}$.
\end{theorem}

\begin{proof} Let us choose an index $j$. Let $F^p\in U(\n)[k^{(j)}]$, and
let $\tilde{F}^p\in U[k^{(j)}]$ be a preimage of $F^p$
under the map $q:U_r\to U(\n)$. First we show that
\begin{description}
\item[(iii)]$<\tilde{F}^{1_j},\bar{\eta}_{\beta_j}>=1$;
\item[(iv)] $<\tilde{F}^p,\bar{\eta}_{\beta_j}>=0$ if 
$p\not=1_j$.
\end{description} For this let
$\overline{\Omega}^{\g}_{k^{(j)}}$ be obtained from
$\Omega^{\g}_{k^{(j)}}$ by replacing $\eta_{\beta_j}$ with
$\bar{\eta}_{\beta_j}$. Using (ii) and definition (\ref{resdef}) of
$<\ ,\
>$ we have $(<\tilde{F}^p,\ \
>\otimes 1)\overline{\Omega}^{\g}_{k^{(j)}}=
(<\tilde{F}^p,\ \
>\otimes 1){\Omega}^{\g}_{k^{(j)}}$. This latter is equal to
\[(1\otimes q)(<\tilde{F}^p,\ \ >\otimes 1)\Omega_{k^{(j)}}=(1\otimes
q)(1\otimes \tilde{F}^p)=1\otimes F^p,\] which proves (iii) and
(iv).

Conditions (i), (iii), and (iv), together with Theorem
\ref{dualcomp} imply $\bar{\eta}_{\beta_j}=\eta_{\beta_j}$.
\end{proof}

\bigskip
\bigskip

\noindent{\bf Notations.} Root systems will be considered in a space
with an ortho-normal basis $\{e_1,\ldots,e_r\}$. For a sequence of
integers $I=(i_1,i_2,\ldots,i_n)$ let $[f_I]$ denote the following
multiple commutator
$[\ldots[[f_{i_1},f_{i_2}],f_{i_3}]\ldots,f_{i_n}]$. Let
$\str(u_1,u_2,\ldots,u_n)$ denote the following tree
\[z\li{1}u_1\li{2}u_2\li{3}\ldots\li{n}u_n.\]

\begin{theorem} \label{lists} The elements given in the following list are the
$\eta_{\beta}$'s of Theorem \ref{product} for the simple Lie
algebras of type $A,B,C,D$ and the choices given.
\end{theorem}
\begin{description}
\item[{\bf A}${}_{r-1}$] For the positive roots
$e_i-e_j$ ($1\leq i<j\leq r$) and simple roots
$\alpha_i=e_i-e_{i+1}$ ($i=1,\ldots,r-1)$, we make the following
choices.
\begin{itemize}
\item $e_i-e_j<e_{i'}-e_{j'}$ if either $i<i'$ or $i=i'$ and $j<j'$;
\item $F_{e_i-e_j}=[f_{(j-1,j-2,\ldots,i)}].$
\end{itemize}
Then we have
$\eta_{e_i-e_j}=\str(t^{(i)}_1,t^{(i+1)}_1,\ldots,t^{(j-1)}_1)$.
\item[{\bf B}${}_r$] For the positive roots $e_i$ ($1\leq i\leq r$),
$e_i-e_j$, $e_i+e_j$ ($1\leq i<j<r$), and simple roots
$\alpha_i=e_i-e_{i+1}$ ($1\leq i<r$) and $\alpha_r=e_r$, we make the
following choices.
\begin{itemize}
\item Let $\beta$ be one of $e_i$, $e_i-e_j$ or $e_i+e_j$, and let $\beta,$ be one of $e_{i'}$, $e_{i'}-e_{j'}$
or $e_{i'}+e_{j'}$. Then we set $\beta<\beta'$ if $i>i'$. If
$i<j<j'$ we set $e_i+e_j<e_i+e_{j'}<e_i<e_i-e_{j'}<e_i+e_j$.
\item $F_{e_i-e_j}=[f_{(j-1,\ldots,i)}]$,
$F_{e_i}=[f_{(r,\ldots,i)}]$,
$F_{e_i+e_j}=(-1)^{\binom{j-i+1}{2}+(r-j+1)}[f_{(i,\ldots,r,r,\ldots,j)}]$.
\end{itemize}
Then we have
$\eta_{e_i-e_j}=\str(t^{(j-1)}_1,t^{(j-2)}_1,\ldots,t^{(i)}_1)$,
$\eta_{e_i}=\str(t^{(r)}_1, t^{(r-1)}_1,\ldots, t^{(i)}_1)$ and
$\eta_{e_i+e_j}=\asym\big( \str(t^{(j)}_2,\ldots, t^{(r)}_2,
t^{(r)}_1,\ldots,t^{(i)}_1)\big)$.
\item[{\bf C}${}_r$] For the positive roots $e_i-e_j$, $e_i+e_j$ ($1\leq i<j\leq r$) and $2e_i$ ($1\leq i\leq r$)
and simple roots $\alpha_i=e_i-e_{i+1}$ ($1\leq i< r$) and
$\alpha_r=2e_r$, we make the following choices.
\begin{itemize}
\item Let $\beta$ be one of $e_i-e_j$, $e_i+e_j$, or $2e_i$, and let
$\beta'$ be one of $e_{i'}-e_{j'}$, $e_{i'}+e_{j'}$, or $2e_{i'}$.
Then we set $\beta<\beta'$ if $i>i'$. For $i<j<j'$ we set
$e_i+e_j<e_{i}+e_{j'}<2e_i<e_{i}-e_{j'}<e_i-e_j$.
\item $F_{e_i-e_j}=[f_{(j-1,j-2,\ldots,i)}]$, $F_{e_i+e_j}=(-1)^{r-j}[f_{(i,i+1,\ldots,r,r-1,\ldots,j)}]$,
$F_{2e_i}=(-1)^{r-i}/2\cdot[f_{(i,i+1,\ldots,r,r-1,\ldots,i)}]$.
\end{itemize}
Then we have $\eta_{e_i+e_j}=\str(t^{(j-1)}_1,\ldots,t^{(i)}_1)$,
$\eta_{e_i+e_j}=\asym\big(\str(t^{(j)}_2,\ldots,t^{(r-1)}_2,$
$t^{(r)}_1,\ldots,t^{(i)}_1)\big)$, and
\[\begin{aligned}
\eta_{2e_i}= &\asym\Big( \begin{picture}(210,25)(0,10)
\put(0,10){$z$} \put(8,10){\li{1}} \put(40,10){$t^{(r)}_1$}
\put(55,15){\line(1,1){15}} \put(60,27){${}_2$}
\put(55,9){\line(1,-1){15}}\put(32,-5){${}_{2r-2i+1}$}
\put(75,32){$t^{(r-1)}_1$}\put(75,-10){$t^{(r-1)}_2$}
\put(105,35){\line(1,0){15}}\put(110,39){${}_3$}
\put(130,35){$\ldots$}
\put(150,35){\line(1,0){30}}\put(155,39){${}_{r-i+1}$}
\put(185,32){$t^{(i)}_1$}
\put(105,-7){\line(1,0){20}}\put(105,-13){${}_{2r-2i}$}
\put(130,-7){$\ldots$}
\put(150,-7){\line(1,0){30}}\put(155,-13){${}_{r-i+2}$}
\put(185,-10){$t^{(i)}_2$}
\end{picture}\Big).
\end{aligned}
\]

\vskip 1 true cm
\item[{\bf D}${}_r$] For the positive roots $e_j-e_i$, $e_j+e_i$ ($1\leq i<j\leq
r)$, and simple roots $\alpha_1=e_1+e_2$, $\alpha_j=e_j-e_{j-1}$
($1<j\leq r$), we make the following choices.
\begin{itemize}
\item Let $\beta$ be one of $e_j-e_i$, $e_j+e_i$, and let $\beta'$
be one of $e_{j'}-e_{i'}$, $e_{j'}+e_{i'}$. Then we set
$\beta<\beta'$ if $j<j'$. For $i'<i<j$ we also set
$e_j+e_i<e_j+e_{i'}<e_j-e_{i'}<e_j-e_i$.
\item $F_{e_j-e_i}=[f_{(j,\ldots,i+1)}]$,
$F_{e_j+e_1}=[f_{(j,\ldots,3,1)}]$,
$F_{e_j+e_i}=[f_{(j,\ldots,1,3,\ldots,i)}]$ ($i>1$).
\end{itemize}
Then we have $\eta_{e_j-e_i}=\str(t^{(i+1)}_1,\ldots,t^{(j)}_1)$,
$\eta_{e_j+e_1}=\str(t^{(1)}_1,t^{(3)}_1,\ldots,t^{(j)}_1)$, and for
$i>1$ we have $\eta_{e_j+e_i}=
\asym\big(\str(t^{(i)}_2,\ldots,t^{(3)}_2,t^{(2)}_1,t^{(1)}_1,t^{(3)}_1,\ldots,t^{(j)}_1)\big)+
\asym\big($ $\str$ $(t^{(i)}_2,\ldots,t^{(3)}_2,$
$t^{(1)}_1,t^{(2)}_1,t^{(3)}_1,\ldots,t^{(j)}_1)\big)$.

\end{description}

\begin{proof} Simple combinatorics shows that for our choices of $F_\beta$ and the order of the $\beta$'s,
Corollary~\ref{coeffs} implies
\[\res^{(i)}_{k^{(j)}_i}\eta_{\beta_j}=\begin{cases}
\eta_{\beta_j-\alpha_i} & \text{if } \beta_j-\alpha_i \text{ is a positive root}, \beta_j-\alpha_i>\beta_j\\
\eta_{{\frac{\beta_j-\alpha_i}{2}}}*\eta_{{\frac{\beta_j-\alpha_i}{2}}}&
\text{if } \frac{\beta_j-\alpha_i}{2} \text{ is a positive root}, \frac{\beta_j-\alpha_i}{2}>\beta_j\\
0 & \text{otherwise.}\end{cases}\] One can check case by case that
the same residue identities hold for the $\eta_{\beta_j}$'s given in
the Theorem, hence condition (ii) of Theorem \ref{barr} is
satisfied.

To prove property (i) of Theorem \ref{barr} recall that the relevant
Serre relations are
\begin{enumerate}
\item \label{1} $[f_{(i,j)}]$ if $(\alpha_i,\alpha_j)=0$,
\item \label{2} $[f_{(i,j,j)}]$ if $(\alpha_i,\alpha_j)=-1$,
\item \label{3} $[f_{(i,j,j,j)}]$ if $(\alpha_i,\alpha_j)=-2$,
$|\alpha_i|<|\alpha_j|$.
\end{enumerate}
If $(\alpha_i,\alpha_j)=0$, then the value of
$(f_If_if_jf_J-f_If_jf_if_J)$ on the $\A(k)$-elements given in the
Theorem are all 0. Indeed, our sign conventions yield that in the
$C_r$ case $f_If_if_jf_J(\eta_{2e_i})=f_If_jf_if_J(\eta_{2e_i})$;
and in all other cases both sides vanish. This proves that on
multiples of (\ref{1}) all $\A(k)$-elements given in the Theorem as
$\eta_{\beta}$'s, vanish. For multiples of (\ref{2}) and (\ref{3})
the proof is similar.
\end{proof}

\section{Differential form representations of $\A(k)$}

A linear map $\phi:\A(k)\to W$ to a vector space $W$ will be called
a representation of $\A(k)$. For a representation $\phi$ we call
$(\phi\otimes 1)\Omega_k$ the canonical $\phi$-element, and
$(\phi\otimes 1)\Omega_k^{\g}=(\phi \otimes q)\Omega_k$ the
canonical $\phi$-element associated with the Lie algebra $\g$. In
this section we consider $\phi$'s for which the canonical
$\phi$-elements associated with simple Lie algebras play an
important role in hypergeometric solutions of KZ-type differential
equations. Theorem \ref{product} together with Theorem~\ref{lists}
give convenient forms of these canonical elements. An example of
application will be shown in Section \ref{scalar}.

Suppose we are given a vector space $W$, and a vector $\phi(T)\in W$
for every ordered $|k|$-tree $T$ on the vertex set $\T(k)\cup\{z\}$.
This data induces a representation $\phi:\A(k)\to W$ if the
assignment $T\mapsto \phi(T)$ respects relations R1 and R2 (see
Section \ref{start}).

\smallskip
\noindent{\bf Notation.} For an edge $e$ of a spanning tree on the
vertex set $\T(k)\cup\{z\}$, $h(e)$ and $t(e)$ denote the head and
tail of the edge $e$, i.e. the vertices adjacent to $e$, farther
resp. closer to $z$.

\subsection{Rational representation.}

Let $W(t^{(i)}_j,z)$ be the vector space of differential forms
in the variables $\T(k)\cup\{z\}$.

\begin{theorem}
Let the edges of an ordered $|k|$-tree be $e_1,\ldots,e_{|k|}$ (in
order). The map
\[T\mapsto\phi(T)=\bigwedge_{i=1}^{|k|} \dlog(h(e_i)-t(e_i))\]
defines a representation $\phi_{rat}:\A(k)\to W(t^{(i)}_j,z)$.
\end{theorem}

\begin{proof}
We need to check the consistency of the definition with relations R1
and R2 from Section~\ref{start}. Relation R1 follows from the
antisymmetry of the $\wedge$-product. In view of R1 it is enough to
check R2 for $a=1$, $b=2$. Let the three distinguished vertices be
$t$, $s$, and $u$. Then we need to check
\[\Big(\frac{d(t-u)}{t-u}\wedge \frac{d(s-t)}{s-t}+
\frac{d(s-u)}{s-u}\wedge \frac{d(t-u)}{t-u}+
\frac{d(t-s)}{t-s}\wedge \frac{d(s-u)}{s-u}\Big)\wedge R=0,\] where
$R$ is the `rest' of the formula (i.e.
$R=\wedge_{i=3}^{|k|}\dlog(h(e_i)-t(e_i))$), but the first factor is
0.


\end{proof}

This representation of the canonical element (as well as the
theorems corresponding to our Theorems \ref{product}, \ref{lists})
is explored in \cite{rsv}. Observe that, according to \cite{arnold},
the rational representation is an isomorphism onto its image. Hence,
the study of the `canonical differential form' of \cite{rsv} is
equivalent to the study of the canonical element of this paper.

\subsection{Theta representation.}\label{theta}

For $z,\tau\in\CC$, Im $\tau>0$, the first Jacobi theta function is
defined by the infinite product
\[\theta(z)=\theta(z,\tau)=ie^{\pi
i(\tau/4-z)}(x;q)(\frac{q}{x};q)(q;q),\qquad q=e^{2\pi i\tau},
x=e^{2\pi i z},\quad (y;q)=\prod_{j=0}^{\infty}(1-yq^j),\]
\cite{ww}. It is an entire holomorphic function of $z$ satisfying
\[\theta(z+1,\tau)=-\theta(z,\tau), \qquad
\theta(z+\tau,\tau)=-e^{-\pi i \tau-2\pi i z}\theta(z,\tau),\qquad
\theta(-z,\tau)=-\theta(z,\tau).\]
By $\theta'(z,\tau)$ we will mean the derivative in the $z$
variable.

\begin{definition}\[\sigma_w(t)=\sigma_w(t,\tau)=\frac{\theta(w-t,\tau)}{\theta(w,\tau)\theta(t,\tau)}\cdot
\theta'(0,\tau).\]
\end{definition}

The listed properties of the theta function yield that the $\sigma$
function---viewed as a function of $t$---has simple poles at the
points of $\Lambda_\tau=\Z\oplus \Z\tau\subset \CC$, as well as the
properties
\begin{equation}\sigma_w(t+1,\tau)=\sigma_w(t,\tau), \qquad
\sigma_w(t+\tau,\tau)=e^{2\pi i w}\sigma(t,\tau),\qquad
\res_{t=0}\sigma_w(t,\tau)=1.\label{sigprop}\end{equation}

\begin{theorem} \label{sigid}
We have
\begin{equation}\label{sigma_identity}
\sigma_{w_1+w_2}(t-u)\sigma_{w_2}(s-t)-\sigma_{w_2}(s-u)\sigma_{w_1}(t-u)+\sigma_{w_1}(t-s)\sigma_{w_1+w_2}(s-u)=0.
\end{equation}
\end{theorem}

\begin{proof} Consider the left hand side of the equation as a
function $f(t)$ of $t$. Then $f$ is a function on $\CC\setminus\big(
(u+\Lambda_\tau)\cup(s+\Lambda_\tau)\big)$, satisfying
$f(t+1)=f(t)$, $f(t+\tau)=e^{2\pi i w_1}f(t)$. The above properties
of $\sigma$ also imply that
\[\begin{aligned}
\res_{t=u} f(t)= & 1\cdot\sigma_{w_2}(s-u)-\sigma_{w_2}(s-u)\cdot
1-0=0 \\
 \res_{t=s} f(t)=&\sigma_{w_1+w_2}(s-u)(-1)+0+1\cdot
\sigma_{w_1+w_2}(s-u)=0.\end{aligned}\] Therefore we have an entire
function. According to $f(t+1)=f(t)$ this can be written in the form
of $\sum d_n e^{2\pi n i t}$. Substituting $f(t+\tau)=e^{2\pi i
w_1}f(t)$ we obtain that $d_n=0$ for all $n$.
\end{proof}

For $k\in \N^r$ we consider variables (weights) $w^{(i)}_j$
($i=1,\ldots,i$, $j=1,\ldots,k_i$), and say that the weight of the
coordinate $t^{(i)}_j\in\T(k)$ is $w^{(i)}_j$.
Let $T$ be a rooted tree (with root $z$) and $v$ a vertex of $T$. We
define the branch $B(v)$ of $v$ to be the collection of those
vertices $w$ for which the unique path connecting $w$ with $z$
contains $v$. By definition $v\in B(v)$.  The load $L(v)$ of a
vertex $v$ in a tree $T$ is defined to be the sum of the weights of
the vertices in $B(v)$. Let $W(t^{(i)}_j,z,w^{(i)}_j,\tau)$ be the
vector space of differential forms in the variables $\T(k)\cup\{z\}$
depending on the parameters $\{w^{(i)}_j\}$ and $\tau$.

\begin{theorem} \label{theta-repr}
Let the edges of an ordered $|k|$-tree be $e_1,\ldots,e_{|k|}$ (in
order). The map
\[T\mapsto \phi(T)=\bigwedge_{i=1}^{|k|}
\sigma_{L(h(e_i))}(h(e_i)-t(e_i))\ d(h(e_i)-t(e_i))\] defines a
representation $\phi_{\theta}:\A(k)\to
W(t^{(i)}_j,z,w^{(i)}_j,\tau)$.
\end{theorem}

\begin{proof}
We need to check the consistency of the definition with relations R1
and R2 from section \ref{start}. Relation R1 follows from the
antisymmetry of the $\wedge$-product.

In view of R1 it is enough to check R2 for $a=1$, $b=2$. In R2, let
the three distinguished vertices be $t$ (bottom), $s$ (upper-right),
$u$ (upper-left). Because of symmetry, we can assume that among
$t,s,u$ it is $u$ that is the closest to the vertex $z$. Consider
the first of the three graphs pictured in R2. Let $L_1$ be the total
weight of the vertices in $B(t)\setminus B(s)$, and let $L_2$ be the
load of $s$, as in the following picture.

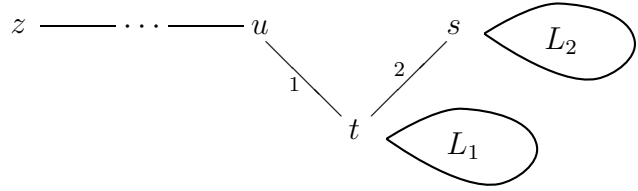
\[
\begin{pspicture}(10,3.5)(0,1)
\put(0,3){$z$}\put(0.4,3.1){\line(1,0){1}} \put(1.5,3){$\cdots$}
\put(2.1,3.1){\line(1,0){1}} \put(3.2,3){$u$}
\put(3.4,2.9){\line(1,-1){1}} \put(4.5,1.6){$t$}
\put(4.8,1.9){\line(1,1){1}} \put(5.8,3){$s$}
\pscurve{-}(5,1.6)(6.5,1)(7,1.5)(6,2)(5,1.6)
\pscurve{-}(6.3,3)(7.8,2.4)(8.3,2.9)(7.3,3.4)(6.3,3)
\put(5.8,1.4){$L_1$} \put(7.1,2.8){$L_2$} \put(3.7,2.3){${}_{1}$}
\put(5.1,2.5){${}_2$}
\end{pspicture}
\]

With these notations we need to check
\[\Big(\sigma_{L_1+L_2}(t-u)d(t-u)\wedge\sigma_{L_2}(s-t)d(s-t)+
\sigma_{L_2}(s-u)d(s-u)\wedge\sigma_{L_1}(t-u)d(t-u)+\]
\[\sigma_{L_1}(t-s)d(t-s)\wedge\sigma_{L_1+L_2}(s-u)d(s-u)\Big)\wedge R=0,\]
where $R$ is the `rest' of the formula, i.e.
$\wedge_{i=3}^{|k|}\sigma_{L(h(e_i))}(h(e_i)-t(e_i))\
d(h(e_i)-t(e_i))$. Considering the path connecting $u$ with $z$ one
finds that $R$ has a $du$ factor, so our equation reduces to
(\ref{sigma_identity}). This proves the theorem.
\end{proof}



\begin{theorem}
The representations $\phi_{rat}$ and  $\phi_\theta$
are injective.
\end{theorem}

\begin{proof} The representation
$\phi_{rat}$ is injective by \cite{arnold}. The representation
$\phi_\theta$ is injective because $\phi_{rat}$ is obtained from
$\phi_\theta$ by the following degeneration. As $\tau \to i\infty$,
the function $\sigma_w(t,\tau)$ tends to $ \pi\sin(\pi(w - t))/(\sin
(\pi w) \sin (\pi t))$. This function tends to $\pi/(e^{2\pi t} -1)$
as $w$ tends to $+\infty$. For small $t$, that function approximates
$2/t$.
\end{proof}

\subsection{The $\g$-theta representation.}\label{gtheta}
Let $\g$ be a simple Lie algebra with Cartan decomposition
$\g=\n\oplus \h \oplus \np$ and simple roots $\alpha_i\in \h^*$.
Consider the theta representation from Section \ref{theta} with the
choices $w^{(i)}_j=\alpha_i(\lambda)$, $\lambda\in \h$. Thus we
obtain a representation $\phi_{\g,\theta}:\A(k)\to
W(t^{(i)}_j,z,\lambda,\tau)$ to the space
$W(t^{(i)}_j,z,\lambda,\tau)$ of differential forms in the variables
$t^{(i)}_j,z$ depending on $\lambda\in \h$ and $\tau\in\CC$, Im
$\tau>0$. We call this the $\g$ theta-representation.

\subsection{The canonical elliptic differential form.}
 For the simple Lie algebra $\g$ we defined a map on
both factors of $\A(k)\otimes U[k]$, namely
$\phi_{\g,\theta}:\A(k)\to W(t^{(i)}_j,z,\lambda,\tau)$ (Section
\ref{gtheta}), and $q:U[k]\to U(\n)[k]$ (Definition \ref{qdef}). Let
the canonical elliptic differential form
$\Theta^{\g}_k=\Theta^{\g}_k(t^{(i)}_j,z,\lambda,\tau)$ be the image
in $W(t^{(i)}_j,z,\lambda,\tau)\otimes U(\n)[k]$ of the canonical
element $\Omega_k$ under the map $\phi_{\g,\theta} \otimes q$.

\begin{example} \rm By Theorems
\ref{product} and \ref{lists} we have
\[
\Theta^{sl_3}_{(2,1)}=(\phi_{sl_3,\theta}\otimes q) \Omega_{(2,1)}=
\]
\[=(\phi_{sl_3,\theta}\otimes 1)\Big(\ \
\asym\Big(
\begin{picture}(65,25)
\put(5,0){$z$}\put(12,7){\line(1,1){10}}\put(13,14){${}_1$}
\put(28,17){$t_1$} \put(36,-15){\line(1,0){15}}\put(42,-10){${}_3$}
\put(53,-20){$s_1$}
\put(12,-3){\line(1,-1){10}}\put(13,-14){${}_2$}\put(28,-20){$t_2$}
\end{picture}\Big)\otimes f_1[f_2,f_1]+\frac{1}{2}\asym\Big(
\begin{picture}(50,25)
\put(5,0){$z$} \put(12,7){\line(1,1){10}}\put(13,14){${}_1$}
\put(28,17){$t_1$}
\put(12,-3){\line(1,-1){10}}\put(13,-14){${}_3$}\put(28,-20){$s_1$}
\put(12,2){\line(1,0){13}} \put(28,0){$t_2$} \put(17,5){${}_{2}$}
\end{picture}
\Big)\otimes f_1^2f_2\ \ \Big)\]\smallskip
\[\begin{aligned} =\Big(\sigma_{\alpha_1(\lambda)}(t_1-z)dt_1\wedge
\sigma_{\alpha_1(\lambda)+\alpha_2(\lambda)}(t_2-z)dt_2\wedge
\sigma_{\alpha_2(\lambda)}(s_1-t_2)ds_1-& & \\
\sigma_{\alpha_1(\lambda)}(t_2-z)dt_2\wedge
\sigma_{\alpha_1(\lambda)+\alpha_2(\lambda)}(t_1-z)dt_1\wedge
\sigma_{\alpha_2(\lambda)}(s_1-t_1)ds_1\Big)& \otimes f_1[f_2,f_1]\\
+\Big(\sigma_{\alpha_1(\lambda)}(t_1-z)dt_1\wedge
\sigma_{\alpha_1(\lambda)}(t_2-z)dt_2\wedge
\sigma_{\alpha_2(\lambda)}(s_1-z)ds_1\Big)& \otimes f_1^2f_2.
\end{aligned}
\]
\end{example}

\section{Application: Eigenfunctions of Calogero-Moser Hamiltonian operator.}\label{scalar}
Let $p$ and $r$ be natural numbers, $\tau\in \CC$, Im $\tau>0$.
Consider the Hamilton operator of the Calogero-Moser quantum
$r+1$-body system
\[H=-\sum_{i=1}^{r+1} \frac{\partial^2}{\partial\lambda_i^2}
- 2p(p+1)\sum_{1\leq i<j\leq
r+1}\rho'(\lambda_i-\lambda_j,\tau),\qquad
\rho(t,\tau)=\frac{\theta'(t,\tau)}{\theta(t,\tau)},\] acting on
scalar functions of $\lambda_1, \ldots, \lambda_{r+1}$, see
\cite{fv}.

\begin{theorem} Let $\alpha_j(\lambda)=\lambda_j-\lambda_{j-1}$ and
$k=(rp,\ldots,2p,p)$. Use the notation $t^{(-1)}_j=0$ for any $j$.
Let $\xi\in \CC^{r+1}$, and suppose that $t^{(i)}_j$ obey the `Bethe
ansatz' equations of \cite[Th.11]{fv}. Then the function
\begin{equation}\label{eigenfunction}
e^{2\pi
i\sum_{j=1}^{r+1}\lambda_j\xi_j}\sym_k\Big(\prod_{l=1}^p\prod_{1\leq
i\leq j\leq r}
\sigma_{(\alpha_i+\ldots+\alpha_j)(\lambda)}(t^{(i)}_{(j-i)p+l}-t^{(i-1)}_{(j-i+1)p+l})
\Big)\end{equation} is an eigenfunction of $H$.
\end{theorem}

\begin{proof} Consider $\g=sl_{r+1}(\CC)$ with its Cartan decomposition and
simple roots as in Theorem \ref{lists}, and the corresponding
$f_i\in \n$. The $p(r+1)$'th symmetric power of the standard
representation of $sl_{r+1}(\CC)$ is the highest weight module
$V_\Lambda$ with highest weight $\Lambda=\sum_i
\frac{(r+1-i)p}{r+1}\alpha_i$.

The surjective map $U(\n)\to V_\Lambda$, $x\mapsto x\cdot v_\Lambda$
induces the map
\[W(t^{(i)}_j,z,\lambda,\tau)\otimes U^{\g}(\n)[k]\to W(t^{(i)}_j,z,\lambda,\tau)\otimes V_\Lambda.\]
Let $\Theta^{V_\Lambda}_k$ be the image of $\Theta^{\g}_k$ under
this map. Theorem 11 of \cite{fv} states that the function
\begin{equation}\label{eigen}
\psi_\gamma(\lambda)=e^{2\pi i\sum_{j=1}^{r+1}  \lambda_j\xi_j}\cdot
\Theta^{V_\Lambda}_k(t,z,\lambda,\tau)\ /\ \bigwedge_{i,j}
dt^{(i)}_j
\end{equation}
is an eigenfunction of the operator $H$, if $z=0$, Im $\tau>0$, and
$t$ satisfies the Bethe ansatz equations. We claim that
\begin{equation}\label{starform}
\begin{aligned}\Theta_k^{V_\Lambda}=
\frac{1}{(p!)^r}\phi_\theta\Big(&\str(z,t^{(1)}_1))^{*p}*
\str(z,t^{(1)}_1,t^{(2)}_1)^{*p}*\ldots\\
& *\str(z,t^{(1)}_1,\ldots,t^{(r)}_1)^{*p}\Big)\otimes
f_1^p[f_2,f_1]^p\ldots
[f_{(r,r-1,\ldots,1)}]^pv_\Lambda.\end{aligned}
\end{equation}
To prove this, choose the ordering of positive roots $\beta$ and the
$F_{\beta}$'s as in Theorem \ref{lists}, and consider
$\Theta_k^{V_\Lambda}$ in the induced PBW form:
\[\sum_p T_p \otimes F_{\beta_1}^{p_1}\cdots F_{\beta_m}^{p_m}v_\Lambda.\]
Observe that the $F_{\beta}v_\Lambda=0$ unless $\beta=e_1-e_j$.
Hence in a non-zero term the powers of the factors $F_{r+1},
F_{r+2}, \ldots$ have to be zero. Therefore there is only one
non-zero term \[\phi_\theta(T_{(p,\ldots,p)})\otimes
f_1^{p}[f_2,f_1]^{p}\ldots[f_{(r,r-1,\ldots,1)}]^{p}v_\Lambda.\]
According to Theorems \ref{product} and \ref{lists} this term is
equal to the right hand side of formula (\ref{starform}).

Tracing back the definitions of the star-product and $\phi_\theta$
in (\ref{starform}), 
we obtain the theorem.
\end{proof}

For $sl_2(\CC)$ the Theorem is proved in \cite{etkir}. A different
Bethe ansatz formula is given in \cite{fv}.

We used the elliptic version of the operator $H$. One can consider
its trigonometric limit as in~\cite{fv2}. To obtain eigenfunctions
in that case we need to replace the representation $\phi_\theta$
with $\lim_{\tau\to i\infty}(\phi_\theta)$ in (\ref{starform}).

\section{Appendix: Notes on the theta representation}

\subsection{The image of the theta representation.}


We give a description of the image of the representation
$\phi_\theta:\A(k)\to W(t^{(i)}_j,z,w^{(i)}_j,\tau)$ for $z=0$ and
$k=(1,\ldots,1)\in \N^r$ and fixed $\tau\in \CC$ with Im $\tau>0$.
The general case can be derived from this special case.

Fix $k=(1,\ldots,1)\in \N^r$, $z=0$, and $\tau\in\CC$ with Im
$\tau>0$. Let $\A(k)_{(z=0)}$ denote the space $\A(k)$ with the
substitution $z=0$. The space $\A(k)_{(z=0)}$ is spanned by ordered
trees on the vertex set $\{t_1,\ldots,t_r,0\}$ (we use the notation
$t_i=t^{(i)}_1$). For such a tree $T$,
\[\phi_\theta(T)=\bigwedge_{i=1}^{r}
\sigma_{L(h(e_i))}(h(e_i)-t(e_i))\ d(h(e_i)-t(e_i))\] is a
meromorphic differential form in the variables $t_1,\ldots,t_r$ and
it depends meromorphically on the complex parameters
$w_1,\ldots,w_r$. Consider the translates of the discriminantal
arrangement $\C^r(0)$ by the points of the lattice
$\Lambda_\tau=\Z+\Z\tau$ in each coordinates. Let $M_r$ be the union
of all these translates, i.e.
\[M_r=\Big(\cup_{H\in
\C^r(0)} H\Big)+\{(x_1,\ldots,x_r): x_j\in\Lambda_\tau\}.\] The
differential form $\phi_\theta(T)$ has at most simple poles at
points of $M_r$ and it is holomorphic on $\CC^r-M_r$.

\begin{theorem} Let $k=(1,\ldots,1)\in\N^r$ and fix $\tau\in \CC$ with Im $\tau>0$.
A meromorphic differential $r$-form $\psi$ in the variables
$t_1,\ldots,t_r$, which depends meromorphically  on the complex
parameters $w_i$, is in the image of $\phi_\theta:\A(k)_{(z=0)}\to
W(t_j,0,w_j,\tau)$ if and only if the following conditions hold.
\begin{enumerate}
\item \label{co1}
$\psi(\ldots,t_j+1,\ldots)=\psi(\ldots,t_j,\ldots)$ for all
$j=1,\ldots,r$.
\item \label{co2} $\psi(\ldots,t_j+\tau,\ldots)=e^{2\pi i w_j}\psi(\ldots,t_j,\ldots)$ for
all $j=1,\ldots,r$.
\item \label{co3} For any subset $\{u_1,\ldots,u_n\}\subset \{t_1,\ldots,t_r\}$, the
differential form
\[\res_{u_1=0}\Big(\res_{u_2=0}\big(\ldots\res_{u_n=0}(\psi)\ldots\big)\Big)\]
has at most simple poles at the points of $M_{r-n}$, and it is
holomorphic on $\CC^{r-n}-M_{r-n}$.
\end{enumerate}
\end{theorem}

\begin{proof} First we show that Properties (\ref{co1})--(\ref{co3})
hold for $\psi\in Im(\phi_{\theta})$.  Property (\ref{co1}) follows
from the first formula in (\ref{sigprop}). To prove
Property~(\ref{co2}) let $T$ be an ordered spanning tree on the
vertex set $\{t_1,\ldots,t_r,0\}$, and let $v_0,\ldots,v_n$ be the
neighbors of $t_j$, among which $v_0$ is the closest to $0$ in~$T$.
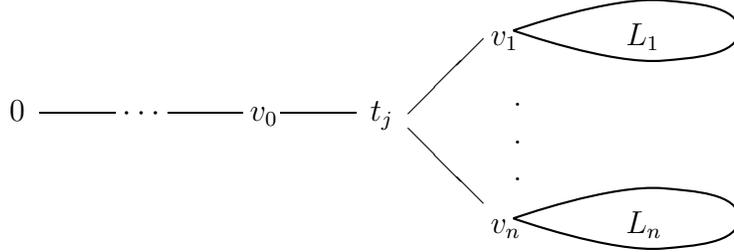
\[
\begin{pspicture}(13,4.4)(0,1.2)
\put(0,3){$0$}\put(0.4,3.1){\line(1,0){1}} \put(1.5,3){$\cdots$}
\put(2.1,3.1){\line(1,0){1}} \put(3.2,3){$v_0$}
\put(3.6,3.1){\line(1,0){1}} \put(4.8,3){$t_j$}
\put(5.3,2.9){\line(1,-1){1}} \put(5.3,3.1){\line(1,1){1}}
\put(6.4,4){$v_1$} \put(6.4,1.5){$v_n$}
\pscurve{-}(6.7,4.2)(8.7,4.6)(9.7,4.2)(8.7,3.8)(6.7,4.2)
\pscurve{-}(6.7,1.7)(8.7,2.1)(9.7,1.7)(8.7,1.3)(6.7,1.7)
\put(8.2,4){$L_1$} \put(8.2,1.5){$L_n$}\put(6.7,2.2){.}
\put(6.7,2.7){.} \put(6.7,3.2){.}
\end{pspicture}
\]
Let $L_i$ be the load of $v_i$. The only factors of $\phi_\theta(T)$
that contain $t_j$ are
\[\sigma_{w_j+\sum
L_l}(t_j-v_0)d(t_j-v_0)\cdot\bigwedge_{l=1}^n
\sigma_{L_l}(v_l-t_j)d(v_l-t_j).\] According to the second formula
from (\ref{sigprop}), after substituting $t^{(i)}_j+\tau$ these
factors get multiplied~by
\[e^{2\pi i(w_j+\sum L_l)} \cdot \prod_{l=1}^n e^{2\pi i
(-L_l)}=e^{2\pi i w_j},\] what we needed to prove. To prove Property
(\ref{co3}) observe that $\res_{t_j=0}
\phi_\theta(T)=\phi_\theta\big( \res^{(j)}_1 T\big).$ Here $\res$ on
the left hand side means residue of differential forms, and on the
right hand side it is the one defined in Section \ref{residue}.
Therefore the multiple residue in Property (\ref{co3}) is the
$\phi_\theta$ image of a certain tree; therefore $\psi$ satisfies
Property (\ref{co3}).
\medskip

We will prove the other direction by induction on $r$. For $k=(1)$
let $\psi(t)$ satisfy the properties of the Theorem. If its residue
at $t=0$ is $c$, then the differential form
\[\psi-c\cdot\sigma_w(t)dt=\psi-\phi_\theta\Big(c\cdot
(0\li{1}t)\Big)\] satisfies Properties (\ref{co1}), (\ref{co2}) and
it is a holomorphic form; thus it has to be 0 (see the argument at
the end of the proof of Theorem \ref{sigid}).

Now let $\psi$ be a differential form in the variables $t_1,\ldots,
t_r$ satisfying the conditions of the theorem. By the induction
hypotheses we have $\res_{t_j=0}\psi =\phi_\theta(T_j)$. Consider
\[\psi-\phi_\theta\Big(\sum_{j} 0\li{1}t_j\li{}T_j  \Big).\]
Here the root vertex $0$ of $T_j$ is glued to the vertex $t_j$ (and
the name $t_j$ is kept for it); the numbers on the edges of $T_j$
are increased by 1. This form has Properties (\ref{co1}),
(\ref{co2}), and its residues at the hyperplanes $t_j=0$ are 0.
Therefore this form is 0 (consider the new variable $\sum t_j$, and
the argument at the end of the Proof of Theorem \ref{sigid}).
Therefore $\psi$ is in the image of $\phi_\theta$.
\end{proof}

\subsection{Differential forms depending on $\tau$.}

We will define a modification of the representation $\phi_\theta$
that takes into account the dependence on the modular parameter
$\tau$ more naturally. Let
\[\omega_w(t)=\sigma_w(t,\tau)dt-\frac{1}{2\pi
i}\partial_w\sigma_w(t,\tau)d\tau.\] This form is used to construct
hypergeometric solutions for elliptic KZ differential equations
in~\cite{fv}.

\begin{theorem}
The differential form $\omega_w(t,\tau)$ is closed in
$(t,\tau)$-space.
\end{theorem}

\begin{proof}
We need to show that \begin{equation}\label{tauident}
(\partial_\tau+\frac{1}{2\pi
i}\partial_t\partial_w)\sigma_w(t,\tau)=0. \end{equation}
Differentiating the functional relation $\sigma_w(t+\tau,\tau)=
\exp(2\pi iw)\sigma_w(t,\tau)$ with respect to $\tau$ and then with
respect to $t,w$, we see that the left-hand side $L(t)$ of
(\ref{tauident}) obeys $L(t+\tau)=\exp(2\pi iw)L(t)$ (the main point
is that the inhomogeneous term cancels). Also trivially
$L(t+1)=L(t)$. The residue at the pole of $\sigma_w$ at $t=0$ is
independent of $w,\tau$ (it is 1) so $L(t)$ is an entire function of
$t$ and must thus vanish.
\end{proof}

\begin{theorem} The form $\omega$ satisfies the following
functional relation
$$
  \omega_{w_1+w_2}(t-u)\wedge\omega_{w_2}(s-t)+
  \omega_{w_2}(s-u)\wedge\omega_{w_1}(t-u)+
  \omega_{w_1}(t-s)\wedge\omega_{w_1+w_2}(s-u)=0.
 $$
\end{theorem}

\begin{proof} The left hand side of the formula is $(f_1dt\wedge
ds+f_2ds\wedge du+f_3du\wedge dt)+(g_1dt+g_2ds+g_3du)\wedge d\tau$.
The first term vanishes by the formula in Theorem \ref{sigid}. The
vanishing of the second term reduces to the vanishing of
$\partial_{w_1}$ and $\partial_{w_2}$ of the formula in Theorem
\ref{sigid}.
\end{proof}

\begin{corollary}
Let the edges of an ordered $|k|$-tree be $e_1,\ldots,e_{|k|}$ (in
order). The map
\[T\mapsto \phi(T)=\bigwedge_{i=1}^{|k|}
\omega_{L(h(e_i))}(h(e_i)-t(e_i))\ d(h(e_i)-t(e_i))\] defines an
injective representation of $\A(k)$ to the space of {\em closed}
differential forms in the variables $\T(k)\cup\{z\}\cup\{\tau\}$,
depending on the parameters $\{w^{(i)}_j\}$.
\end{corollary}

\end{document}